\documentclass{article}
\usepackage{amssymb, amsmath}
\usepackage{systeme}
\usepackage{mathtools}
\DeclarePairedDelimiter{\ceil}{\lceil}{\rceil}

\begin{document}

\title{CANONICAL DECOMPOSITION  IN $\mathbb{R}_{+}^{*}$ OF A CONVERGENT  NATURAL NUMBER BY THE COLLATZ ITERATIONS }
\author{Esse Koudam \\
UNC Charlotte, Charlotte NC 28223, USA}

\maketitle

\begin{abstract}

The Collatz variations pattern seems not to have any recurrence relation between numbers. But knowing that there is at least a natural number that converges after several iterations we construct a function $f_{X,Y}$ that is equal to the value of convergence for all convergent sequences. A canonical decomposition can be expressed for such numbers.

\end{abstract}
$\mathbf{Keywords}$: Syracuse problem, 3N+1 problem, Collatz conjecture. 

\section{Introduction}

The Collatz conjecture has many denominations. It is also known as the Syracuse problem or the 3N+1 problem. The problem was first stated by the German mathematician Lothar Collatz in the 1930's \cite{Lagarias}. The conjecture is sumarized as follows. Take any natural number $n$ not equal to zero. If $n$ is even divide by 2. If $n$ is odd multiply it by 3 and add 1. Repeat the process to infinity. Does the sequence created reaches 1 for every initial number $n$?  
The Collatz sequence $(C_{p})_{p\in \mathbb{N}}$ started with a natural number $n$ different of zero is called convergent when after $k$ iterations the sequence is equal to 1. The total stopping time $\sigma_{\infty}(n)=inf\left\{k: T^{k}(n)=1\right\}$ \cite{Marc Chamberland}; $k$ is the finite least iterations before $(C_{p})$ converges. Condider the function:
\[
g(n) = \begin{cases} \dfrac{n}{2}, & \mbox{if } n\mbox{ is even} \\ 3n+1, & \mbox{if } n\mbox{ is odd} \end{cases}
\]
Form the sequence by performing an infnite operation of the fonction.
Notation:
\[
C_{p} = \begin{cases} n, & \mbox{for } p=0 \\ g\left(C_{p-1}\right), & \mbox{for } p>0 \end{cases}
\]
$C_{p}$ is the value of $g$ applied to $n$ recursively $p$ times: in notation $C_{p}=g^{p}\left(n\right)$. The smallest $p$ such that $C_{p}=1$ is nothing than $\sigma_{\infty}(n)$ defined earlier as the total stopping time ($p=k$). \\
\\
A divergent sequence isn't yet found.  The divergence  would consist of  a total stopping being infinity. In notation: $\sigma_{\infty}(n)=\infty$ \cite{Marc Chamberland}.  Even though computational method had proven the convergence of  all natural number $n< 20\cdot2^{58}$ \cite{Silva}, does not totally prove the Collatz conjecture. But it tells us the existence of several convergent numbers (The partition set of the convergent numbers in $\mathbb{N}$  is not empty). \\
\\
This document is intented to prove that all convergent numbers have their convergence same as a function $f_{X,Y}$ . In general, $g^{k}\left(n\right)=f_{X,Y}=1$. This paper also includes properties of convergent numbers by the Collatz sequence and a generalisation of the idea that the set of convergent  $n$  is never empty to an infinite set.

\section{The odd and even iterations X and Y at convergence}

\subsection{The $k-tuple$ associated to a Collatz sequence at the total stopping time $\sigma_{\infty}(n)=k$ }

$\mathbf{Definition\mbox{ }2.1}$: The $k-tuple$ associated to $n$ after $k$ total Collatz  iterations is the chain constitued of all values $C_{p}$ when $p$ varies from $0$ to $p-1$ ($p=0,1,......,k-1$). Then the $k-tuple$ is $(C_{0},C_{1},.......,C_{k-1})$. \\
\\
$\mathbf{Examples\mbox{ }2.1}$ \\
For $n=6$, $C_{8}=1$. The $8-tuple$ associated to $6$ is $(6, 3, 10, 5, 16, 8, 4, 2)$ . \\
\\
Another example is $n=19$, it takes $20$ iterations before it gets to $1$. $C_{20}=1$ and its $20-tuple$ is $(19, 58, 29, 88, 44, 22, 11, 34, 17, 52, 26, 13, 40, 20, 10, 5, 16, 8, 4, 2)$ . \\

\subsection{The smallest odd and even iterations $X$ and $Y$ }
Consider the $k-tuple$ $(C_{0},C_{1},.......,C_{k-1})$ associated to a convergent sequence of $n$. Let us make a set $E$ of all the element in the $k-tuple$ chain and 2 subsets $E_{1}$, $E_{2}$ defined respectively as set of all the odd and all even numbers of $E$. \\
\begin{equation}
E = \left\{  {C_{0},C_{1},.......,C_{k-1}} \right\} .
\end{equation}
$\mathbf{Definition\mbox{ }2.2}$: The iterations on odd numbers X is the cardinal of the set $E_{1}$ and $Y$, the  iterations on even number is the cardinal of the set $ E_{2}$.\\
\begin{equation}
 Card\left\{ E_{1}\right\}=X
\end{equation}
and,
\begin{equation}
 Card\left\{ E_{2}\right\}=Y .
\end{equation}
\\
$\mathbf{Remark\mbox{ }2.2}$ \\
$Card\left\{ {E} \right\} = Card\left\{ E_{1}\right\}+ Card\left\{ E_{2} \right\}$, \\
$\left\{ {E} \right\}= E_{1}\cup E_{2}$, \\
$\sigma_{\infty}(n)=k=X+Y$ . \\
By convenience we'll note a convergent sequence of $n$ after $X$ and $Y$ iterations $n=n(X,Y)$, and we'll denote by $\mathcal{N}$ the set of convergent natural numbers. \\
By definition of the Collatz conjecture $n\not= 0$, so $\mathcal{N} \subseteq \mathbb{N}^{*}$.

\section{The function $f_{X,Y}$ associated to the convergence value $C_{k}(n)$ of $n$}
\subsection{The value at the convergence}
$\mathbf{Definition\mbox{ }3.1}$: The sequence $(C_{p})_{p \in \mathbb{N}}$ is called convergent when after $p= k$ iterations $C_{k}(n)=1$. The value at the convergence of any Collatz sequence started with $n$ non-zero positive integer is the limit taken at the total stopping time. In terms of limit notation:
\medskip
\[ \lim_{p \to \sigma_{\infty}} (C_{p}) = 1 .\]

\subsection{The function $f_{X,Y}$}
$\mathbf{Definition\mbox{ }3.2}$: Let $\mathcal{N}$ be the set  ol the convergent $n=n(X,Y)$, with the couple $(X,Y)$ associated to the its respective $n$. By a function at the convergence of $n$ we mean a map
\[
f_{X,Y}: \mathcal{N} \to \mathbb{N}
\]
where:
\begin{equation}
f_{X,Y}\left(n\right)=\ceil[\big]{\dfrac {3^{X}\left(2n+1\right)-1}{2^{Y+1}}} .
\end{equation}
\\
$\mathbf{Lemma\mbox{ }3.1}$: Let $Z_{n}= \dfrac {3^{X}\left(2n+1\right)-1}{2^{Y+1}}$.  $\forall{n=n(X,Y)} \in \mathcal{N}$, $\exists \varepsilon$ such that $0 \leq \varepsilon <\dfrac{1}{3}$ , for which,
\begin{equation}
 Z_{n}=1-\varepsilon .
\label{eq:five}
\end{equation}
$Proof$: For $n(0,i)=2^{i}$ ($i \in \mathbb{Z_{+}^{*}}$), then $Z_{n}=\dfrac {2^{i}}{2^{i}}=1$ where $\varepsilon =0$. \\
Let $n\not=n(0,i)$. Proceed by ABSURD i.e we supose that $\exists n=n(X,Y)$, and $\varepsilon>\dfrac{1}{3}$ such that $Z_{n}\not=1-\varepsilon$ \\
If $Z_{n}<1-\varepsilon$ \\
$Z_{n}<1-\varepsilon \Longrightarrow \exists \varepsilon'$ such that $\dfrac{1}{3}<\varepsilon'<1$ and $Z_{n}=1-\varepsilon'$. That's ABSURD. \\
If $Z_{n}>1-\varepsilon$ \\
$Z_{n}>1-\varepsilon \Longrightarrow \exists \varepsilon'$ such that $0<\varepsilon'<\varepsilon<\dfrac{1}{3}$ and $Z_{n}=1-\varepsilon'$ (ABSURD), \\
or \\
$\exists L>1$ and $\varepsilon''$, for which $Z_{n}=L-\varepsilon''$ with $\dfrac{1}{3}< \varepsilon''<1$ ($L \in \mathbb{N}$).
\begin{align*}
3Z_{n}+1 &=3L-3\varepsilon'' +1 \\
                 &=(3L+1)-3\varepsilon'' \\
\end{align*}
$n$ converges, $n=n(X,Y) \Longrightarrow  Z_{n}=1-\varepsilon''$ and $0< \varepsilon'' <\dfrac{1}{3}$ ($n\not=n(0,i)$) .
\begin{align*}
3Z_{n}+1 &=3-3\varepsilon'' +1 \\
                 &=4-3\varepsilon'' \\
\end{align*}
For the same $n=n(X,Y)$, we have 2 values of $3Z_{n}+1$ where one's function of $L$. Since a number cannot differ from itself, $L$ must be equal to $1$ and $0<\varepsilon''<\dfrac{1}{3}$. There is a contradiction meaning that there is no such $ L$ greater than  $1$ and there is no $\varepsilon''$, such that $\dfrac{1}{3}< \varepsilon''<1$ for which $Z_{n}=L-\varepsilon''$ . \\
\\
$Conclusion$: $\forall{n=n(X,Y)} \in \mathcal{N}$, $\exists \varepsilon$ such that $0 \leq \varepsilon <\dfrac{1}{3}$ for which $Z_{n}=1-\varepsilon$.\\
\\
We can now prove the following theorem; \\
\\
$\mathbf{Theorem\mbox{ }3.1}$: For all $n=n(X,Y) \in \mathcal{N}$, there is a function $f_{X,Y}(n)=\ceil[\big]{\dfrac {3^{X}\left(2n+1\right)-1}{2^{Y+1}}}$ which is equal to the value of convergence $C_{k}$ of $n$. \\
\begin{equation}
\forall{n=n(X,Y)} \in \mathcal{N}, f_{X,Y}(n)=1 . 
\label{eq:six}
\end{equation}
$Proof$: From Lemma 3.1 for all $n=n(X,Y)$ there is always $\varepsilon$ satisfying the condition $0 \leq \varepsilon <\dfrac{1}{3}$ and we have $Z_{n}=1-\varepsilon$. \\
\begin{align*}
Z_{n}=1-\varepsilon                              &\Longrightarrow \ceil[\big]{Z_{n}}=1 , \\
f_{X,Y}\left(n\right)=\ceil[\big]{Z_{n}} &\Longrightarrow f_{X,Y}\left(n\right)=1 . \\
\end{align*} 
\\
$\mathbf{Corollary \mbox{ }3.1}$: $\forall{n=n(X,Y)} \in \mathcal{N}$, $f_{X,Y}(n)$  is a constant function.

\section{Canonical decomposition $n'$ of $n=n(X,Y)$}
$\mathbf{Definition\mbox{ }4.1}$: For $n$ in $\mathcal{N}$ with its respective couple $(X,Y)$, the expression $n'$  in $\mathbb{R}_{+}^{*}$ of $n$ is:
\begin{equation}
n'=2^{Y}(3^{-X}-\varepsilon_{n}); \mbox{ }
\end{equation}  
$n'$ is called by definition the canocical decomposition of $n$. \\
$\mathbf{Lemma\mbox{ }4.1}$: For all convergent $n=n(X,Y)$, $0<\dfrac{1}{2}(1-3^{-X})< \dfrac{1}{2}$ . \\
$Proof$: $\forall X\in \mathbb{N}$,
\begin{align*}
0 &<3^{-X} <1 \\
-1&<-3^{-X} <0 \\
0 &<1-3^{-X} <1 \\
0 &<\dfrac{1}{2}(1-3^{-X}) < \dfrac{1}{2} \\
\end{align*} \\
$\mathbf{Lemma\mbox{ }4.2}$: If $\left\{n'\right\}$ is the fractionnal part of $n'$ then, $\left\{n'\right\}= \dfrac{1}{2}(1-3^{-X})$ and \mbox{ } $\varepsilon_{n}= 3^{-X}\varepsilon$ . \\
$Proof$: Let $n(X,Y) \in \mathcal{N}$, then from (\ref{eq:six}) we know that $\ceil[\big]{Z_{n}}=1$ .
\[
\ceil[\big]{Z_{n}}=1 \Longleftrightarrow Z_{n}=1-\varepsilon
\]
\begin{align*}
\dfrac{3^{X}\left(2n+1\right)-1}{2^{Y+1}} &=1-\varepsilon \\
                               3^{X}\left(2n+1\right)-1 &=2^{Y+1}(1-\varepsilon) \\
                                                                        &=2^{Y+1}-2^{Y+1}\varepsilon \\
                                 2\cdot3^{X}n+3^{X}-1 &=2^{Y+1}-2^{Y+1}\varepsilon \\
                                                2\cdot 3^{X}n &=1-3^{X}+2^{Y+1}-2^{Y+1}\varepsilon \\
                                                                     n &=\dfrac{1}{2}(3^{-X})-\dfrac{1}{2}+3^{-X}\cdot2^{Y}-(3^{-X}\varepsilon)\cdot2^{Y} \\
                                                                        &=\dfrac{1}{2}(3^{-X}-1)+3^{-X}\cdot2^{Y}-2^{Y}\varepsilon_{n} \\
                                                                        &=\dfrac{1}{2}(3^{-X}-1)+ 2^{Y}(3^{-X}-\varepsilon_{n}) \\
                                                                        &=2^{Y}(3^{-X}-\varepsilon_{n})-\dfrac{1}{2}(1-3^{-X}) \\
                                                                        &= n'- \dfrac{1}{2}(1-3^{-X})\\
\end{align*}
\begin{equation}
n= n'- \dfrac{1}{2}(1-3^{-X}) .
\label{eq:eight}
\end{equation}
The number $n$ is a natural number written as the difference of $2$ real numbers which are positive. Also $0<\dfrac{1}{2}(1-3^{-X})<\dfrac{1}{2}$ and $n'> \dfrac{1}{2}(1-3^{-X})$. The relation $n= n'- \dfrac{1}{2}(1-3^{-X})$ is true if and only if $\dfrac{1}{2}(1-3^{-X})$ is the fractionnal part of $n'$; i.e $\left\{n'\right\}= \dfrac{1}{2}(1-3^{-X})$. \\
$\mathbf{Lemma\mbox{ }4.3}$: If $n'$ is the canocical decomposition of $n$ then $0 \leq \varepsilon_{n} <\dfrac{1}{3^{X+1}}$  . \\
$Proof$: $\varepsilon_{n}= 3^{-X}\varepsilon$ and $0 \leq \varepsilon <\dfrac{1}{3}$ . \\
$\mathbf{Remarks\mbox{ }4.1}$:
For the same iterations $X$ and $Y$ at the convergence of $n$ and $m$, $\varepsilon_{n} =\varepsilon_{m}$ , \\
For $n=n(0,i)$ (or $n=2^{i}$, $i \in \mathbb{Z_{+}^{*}}$), $\varepsilon_{n}=0$ . \\
\\
We can now state the following theorem; \\
\\
$\mathbf{Theorem\mbox{ }4.1}$: Let $n'$ be the canonical decomposition of $n=n(X,Y)$ in $\mathbb{R}_{+}^{*}$. The expression of $n$ in function of $n'$ is:
\begin{equation}
n=\lfloor{n'}\rfloor .
\end{equation}
$Proof$: From (\ref{eq:eight}) we have the equality $n= n'- \dfrac{1}{2}(1-3^{-X})$
\begin{align*}
n                        &= n'- \dfrac{1}{2}(1-3^{-X}) \\
n'                       &= n+ \dfrac{1}{2}(1-3^{-X}) \\
\lfloor{n'}\rfloor&= \lfloor{n+ \dfrac{1}{2}(1-3^{-X})}\rfloor \\
\lfloor{n'}\rfloor&=\lfloor{n+ \left\{n'\right\}}\rfloor \\
\lfloor{n'}\rfloor&= \lfloor{n}\rfloor \\
\lfloor{n'}\rfloor&= n . \\
\end{align*}
\\
$\mathbf{Corollary \mbox{ }4.1}$: $\forall{n} \in \mathcal{N}$, $n=n(X,Y) \Longleftrightarrow n=\lfloor{n'}\rfloor$ . \\

\subsection{Properties}
Consider $\mathcal{N}$, $n'$ the canocical decomposition of $n=n(X,Y)$, and $a$ and $b$ be 2 elements of $\mathcal{N}$. We consider the following strong properties arising from the canonical decomposition: \\
\\
$\mathbf{Unicity \mbox{ } of \mbox{ } the \mbox{ } couple \mbox{ } (X,Y)}$: $a=a(X,Y)$ and $b=b(X,Y)$ iff $a=b$ .\\
\\
$Proof$: Let $a=a(X,Y)$ and $b=b(X,Y)$ then $a=\lfloor{a'}\rfloor$ and $b=\lfloor{b'}\rfloor$
\begin{align*}
a-b &=\lfloor{a'}\rfloor-\lfloor{b'}\rfloor , \\
a-b &= \lfloor2^{Y}(3^{-X}-\varepsilon_{a})\rfloor - \lfloor2^{Y}(3^{-X}-\varepsilon_{b})\rfloor ,
\end{align*}
It's known from remark 4.1 that for the same iterations $X$ and $Y$ at the convergence of $n$ and $m$, $\varepsilon_{n} =\varepsilon_{m}$ then:
\begin{align*}
\varepsilon_{a} &=\varepsilon_{b} , \\
                    a-b &=0 , \\
                       a &=b . \\
\end{align*}
$\mathbf{The \mbox{ } a+b\mbox{ } addition}$: If $a$ and $b$ converge so does $a+b$: i.e $a=a(X,Y)$ and $b=b(X',Y')$, then $\exists (X'',Y'')$ such that $a+b=[a+b](X'',Y'')$ \\
\\
$Proof$: Let $a=a(X,Y)$ and $b=b(X,Y)$ then $a=\lfloor{a'}\rfloor$ and $b=\lfloor{b'}\rfloor$
\begin{align*}
a+b &=\lfloor{a'}\rfloor+\lfloor{b'}\rfloor , \\
a+b &=\lfloor{a'+b'}\rfloor .
\end{align*}
Because $\left\{a'\right\}$ and $\left\{b'\right\}$ are both less than $\dfrac{1}{2}$ . (See Lemma 4.1). \\

\section{Algebra of the set $\mathcal{N}$}
\subsection{Equipotence to $\mathbb{N}$}
$\mathbf{Definition \mbox{ } 5.1}$: $\mathcal{N}$ is equipotent to $\mathbb{N}$ or countably infinite when there exist a function bijective from $\mathcal{N}$ to $\mathbb{N}$. \\
\\
$\mathbf{Lemma \mbox{ } {5.1}}$: There is a bijection from $\mathbb{N}$ to its infinite subsets especially to $\mathcal{N}$. \\
\\
$Proof$: $\mathcal{N} \subset \mathbb{N}$ ; Let consider an order relation $\leq$ on $\mathcal{N}$ and let the set be finite,
\[
\exists M\in \mathcal{N} | M=\left\{k | \forall{x}\in  \mathcal{N}, x< k\right\}
\]
\[
M+1 \notin \mathcal{N}
\]
From properties above the addition of 2 convergent numbers is convergent so must $M+1$  also be in $\mathcal{N}$ i.e also convergent. We arrive at a contracdition. $\mathcal{N}$ is not majored and not a finite set. \\
The application which to every single element of $\mathcal{N}$ associate their perfect square in $\mathbb{N}$  is bijective.

\subsection{Total order relation in $\mathcal{N}$} 
The order relation $\leq$ is total in $\mathcal{N}$. By definition $\leq$ is a total relation order when $\forall a$ and $b$ in the set such that  $a\leq b$, there is also $c$ in the set such $a+c=b$.\\
In $\mathcal{N}$ this relation is verified. In fact, if there is $M$ in $\mathcal{N}$, $M+1$ also is in $\mathcal{N}$ leading to state that two elements $a$ and $b$ in the set are always comparable: $a\leq b$ or $b\leq a$. \\
\\
\subsection{Conclusion}
$We\mbox{ } recall\mbox{ } some\mbox{ } basic\mbox{ }properties$:
Any partition of $\mathbb{N}$ different from the empty set ($\emptyset$) has a least element. The least element to converge in $\mathcal{N}$ is $1$ . \\
The addition in $\mathcal{N}$ is an internal law of composition. \\
So we can assure these following inclusions:
\begin{equation}
\mathcal{N} \subset \mathbb{N}^{*}, and\mbox{ }\\
\mathbb{N}^{*} \subset \mathcal{N} \\
\end{equation}

\end{document}